\documentclass[12pt]{amsart}
\usepackage{amssymb}
\input amssym.def
\usepackage{amsmath, amsfonts,hyperref}
\usepackage{amscd}
\usepackage[mathscr]{eucal}
\newcommand{\Z}{{\mathbb Z}}

\newtheorem{thm}{Theorem}
\newtheorem{lem}{Lemma}
\newtheorem{cor}{Corollary}
\newtheorem{prop}{Proposition}
\newtheorem{rmk}{Remark}
\newtheorem{defn}{Definition}
\newcommand{\thmref}[1]{Theorem~\ref{#1}}
\newcommand{\propref}[1]{Proposition~\ref{#1}}
\newcommand{\lemref}[1]{Lemma~\ref{#1}}
\newcommand{\corref}[1]{Corollary~\ref{#1}}

\begin{document}

\title[Error term in a Parseval type formula II]
{On the error term in a Parseval type formula in the theory of
Ramanujan expansions II}

\author{Giovanni Coppola, M. Ram Murty and Biswajyoti Saha}

\address{Giovanni Coppola\\ \newline
Universit\'a degli Studi di Napoli,
Complesso di Monte S. Angelo-Via Cinthia
80126 Napoli (NA), Italy}
\email{giovanni.coppola@unina.it}
\thanks{Research of the first author was supported by an INdAM Research Grant
(titolare di un Assegno ``Ing. Giorgio Schirillo" dell'Istituto
Nazionale di Alta Matematica)}

\address{M. Ram Murty\\ \newline
Department of Mathematics, Queen's University,
Kingston, Ontario, K7L 3N6, Canada}
\email{murty@mast.queensu.ca}
\thanks{Research of the second author was partially supported by an
NSERC Discovery grant.}

\address{Biswajyoti Saha\\ \newline
Institute of Mathematical Sciences, C.I.T. Campus, Taramani, 
Chennai, 600 113, India}
\email{biswajyoti@imsc.res.in}



\begin{abstract}
For two arithmetical functions $f$ and $g$ with absolutely convergent Ramanujan
expansions, Murty and Saha have recently derived asymptotic formulas with error
term for the convolution sum $\sum_{n \le N} f(n) g(n+h)$ under some suitable
conditions. In this follow up article we improve these results with a weakened
hypothesis which is in some sense minimal.
\end{abstract}

\maketitle

\section{Introduction}
In the seminal article \cite{SR}, Ramanujan unfolded the theory of
Ramanujan sums and Ramanujan expansions. He defined:

\begin{defn}
For positive integers $r,n$,
$$
 c_r(n):=\sum_{a\in (\Z/r\Z)^*}\zeta_r^{an},
$$
where $\zeta_r$ denotes a primitive $r$-th root of unity.
\end{defn}

Since then, these sums are attributed to Ramanujan and
called Ramanujan sums. It is not hard to write $c_r(n)$ in
terms of the M\"obius function $\mu$ (see \cite{RM}). One has
\begin{equation}\label{crn-exp1}
c_r(n)=\sum_{d|n, d|r} \mu(r/d) d.
\end{equation}
One also has the following explicit formula due to H\"older:
\begin{equation}\label{crn-exp2}
c_r(n)=\frac{\varphi(r)}{\varphi(r/d)} \mu(r/d),
\end{equation}
where $d = \gcd(n,r)$ and $\varphi(\cdot)$ denotes the Euler's totient function.

Ramanujan studied these sums in the context of
point-wise convergent series expansion of the form
$\sum_r a_rc_r(n)$ for various arithmetical functions.
Such expansions are now known as Ramanujan expansions. More precisely:

\begin{defn}
We say an arithmetical function $f$
admits a Ramanujan expansion, if for each $n$,
$f(n)$ can be written as a convergent series
of the form
$$
f(n)= \sum_{r \ge 1} \hat f(r)c_r(n)
$$
for appropriate complex numbers $\hat f(r)$.
The number $\hat f(r)$ is said to be the $r$-th Ramanujan coefficient
of $f$ with respect to this expansion.
\end{defn}

Ramanujan himself observed \cite{SR} that such an expansion
is not necessarily unique. He remarked that the assertion
$$
 \sum_{r \ge 1} \frac{c_r(n)}{r} = 0
 $$
is equivalent to the prime number theorem. This equation can be
viewed as a Ramanujan expansion of the zero function.

The vast archive on the theory of Ramanujan sums and Ramanujan expansions
symbolize its great developments, in many directions with many different
aspects, over the past 100 years or so. 
For instance, shortly after Ramanujan's death, Hardy \cite{hardy} proved that
$$ {\phi(n) \over n} \Lambda(n) = \sum_{r=1}^\infty {\mu(r)\over \phi(r)}c_r(n), $$
where $\Lambda$ is the von Mangoldt function.  The series on the righ hand
side is conditionally convergent and so is difficult to use.  
In \cite{padma}, the authors showed that if we ignore convergence questions,
Hardy's formula can be used to derive the Hardy-Littlewood conjecture
about prime tuples. More precisely, one can derive the heuristic
result that
$$\sum_{n\leq N} \Lambda(n) \Lambda(n+h) \sim N \sum_{r=1}^\infty  {\mu^2(r)
\over \phi(r)} c_r(h), $$
a conjecture formulated by Hardy and Littlewood using the more difficult
circle method of Ramanujan.  This led the authors of \cite{GMP} to study
convolution
sums of the kind $\sum_{n \le N} f (n)g(n + h)$ for two arithmetical
functions $f$ and $g$ with absolutely convergent Ramanujan expansions.
They derived asymptotic formulas for such sums, which are analogous to
Parseval's formula in the case of Fourier series expansions. However, the
study of the error term for such formulas was not carried out there. It was
then addressed in \cite{MS} by Murty and Saha. Under certain extended
hypotheses they provide explicit error terms for such formulas.
\par
The works \cite{GMP} and \cite{MS} had some severe
restrictions on the growth of the Ramanujan coefficients.  The goal of 
this paper is to relax these conditions.
To be precise, in \cite{GMP} the authors had the following condition on
the Ramanujan coefficients of $f$ and $g$:
$$
\sum_{r,s} \big| \hat f(r) \hat g(s)\big|(rs)^{1/2}d(r)d(s) < \infty,
$$
which was then extended as
$$
\big|\hat f(r)\big|,\big|\hat g(r)\big| \ll \frac{1}{r^{1 +\delta}}
~\mbox{ for }~ \delta > 1/2,
$$
in \cite{MS}. The condition `$\delta > 1/2$', a priori seems to be ad hoc,
as a similar condition for $\delta > 0$ would be sufficient to ensure the
absolute convergence of the Ramanujan expansion. But if one wants to
extend the hypothesis of \cite{GMP}, in the form that they have in \cite{MS},
then `$\delta > 1/2$' is somewhat an optimal choice.

More recently, Saha \cite{BS} has considered the single sum
$\sum_{n \le N} f(n)$ for an arithmetical function $f$ with absolutely
convergent Ramanujan expansion, in the context of deriving an asymptotic
formula with explicit error term for such a sum. The author obtains his
result under the above condition with $\delta >0$. He also exhibits that
even with the stronger condition that $\delta > 1/2$, one may end up
getting a weaker result if the concerned sum is not handled carefully
enough. This is exactly the phenomenon that we address here.

Hence, for our purpose we enforce the weakened (essentially a minimal)
hypothesis `$\delta>0$' (compared to both \cite{GMP} and \cite{MS}),
and still obtain a better error term. In the following section we state our
results and compare them with their predecessors.

\section{Statements of the Theorems}

In this article we prove the following theorems:

\begin{thm}\label{error1}
 Suppose that $f$ and $g$ are two arithmetical functions with absolutely
 convergent Ramanujan expansions:
 $$
 f (n) = \sum_r \hat f(r) c_r(n), \phantom{mm} g(n) = \sum_s \hat g(s)c_s(n),
 $$
 respectively. Further suppose that
 $$
 \big|\hat f(r)\big|,\big|\hat g(r)\big| \ll \frac{1}{r^{1 +\delta}}
 $$
 for some $\delta > 0$. Then for a positive integer $N$, we have,
 $$
 \sum_{n \le N} f(n) g(n) =
 \begin{cases}
N \sum_{r \ge 1} \hat f(r) \hat g(r) \varphi(r) +
O(N^{1-\delta}(\log N)^{4-2\delta}) & ~ \mbox{ if } ~ \delta < 1,\\
N \sum_{r \ge 1} \hat f(r) \hat g(r) \varphi(r) +
O(\log^3 N) & ~ \mbox{ if } ~ \delta = 1,\\
N \sum_{r \ge 1} \hat f(r) \hat g(r) \varphi(r) +
O(1) & ~ \mbox{ if } ~ \delta > 1.
 \end{cases}
 $$
\end{thm}

\begin{rmk}
\rm The exponent of $N$ in the error term in the above theorem
is consistent with error term which has been obtained in \cite{BS} for sum of
the form $\sum_{n \le N} f(n)$, which was not the case for results in \cite{MS}.
\end{rmk}

\begin{thm}\label{error2}
 Let $f$ and $g$ be two arithmetical functions with the same hypotheses as
 in \thmref{error1} and $h$ be a positive integer.
 Then we have,
  $$
 \sum_{n \le N} f(n) g(n+h) =
 \begin{cases}
N \sum_{r \ge 1} \hat f(r) \hat g(r) c_r(h) +
O(N^{1-\delta}(\log N)^{4-2\delta}) & ~ \mbox{ if } ~ \delta < 1,\\
N \sum_{r \ge 1} \hat f(r) \hat g(r) c_r(h) +
O(\log^3 N) & ~ \mbox{ if } ~ \delta = 1,\\
N \sum_{r \ge 1} \hat f(r) \hat g(r) c_r(h) +
O(1) & ~ \mbox{ if } ~ \delta > 1.
 \end{cases}
 $$
\end{thm}

\thmref{error1} and \thmref{error2} are substantially improved versions,
in terms of both hypotheses and the conclusion, of the following two
theorems respectively.

\begin{thm}[Murty-Saha]\label{error3}
 Let $f$ and $g$ be two arithmetical functions as in \thmref{error1},
 with the last assumption being replaced by the condition that
 $$
 \big|\hat f(r)\big|,\big|\hat g(r)\big| \ll \frac{1}{r^{1 +\delta}}
 $$
 for some $\delta > 1/2$. Then one has,
 $$
 \sum_{n \le N} f(n)g(n) = N \sum_r \hat f(r) \hat g(r) \varphi(r) +
 O\left( N^{\frac{2}{1+2\delta}} (\log N)^{\frac{5+2\delta}{1+2\delta}}\right).
 $$
\end{thm}

\begin{thm}[Murty-Saha]\label{error4}
 Let $f$ and $g$ be two arithmetical functions with the same hypotheses as
 in \thmref{error3} and $h$ be a positive integer. Then,
 $$
 \sum_{n \le N} f(n)g(n+h) = N \sum_r \hat f(r) \hat g(r) c_r(h) +
 O\left( N^{\frac{2}{1+2\delta}} (\log N)^{\frac{5+2\delta}{1+2\delta}}\right).
 $$
\end{thm}

By virtue of \thmref{error2}, we can now naturally extend and improve
the corollaries that were obtained in \cite{MS}. We quote:

\begin{cor}\label{c1}
 For $s,t > 0$, let $\delta:=\min \{s,t\}$. Then for any positive
 integer $h$, we have,
$$
 \sum_{n \le N} \frac{\sigma_s(n)}{n^s}~\frac{\sigma_t(n+h)}{(n+h)^t}=
 \begin{cases}
N\frac{\zeta(s+1)\zeta(t+1)}{\zeta(s+t+2)} \sigma_{-(s+t+1)}(h) +
O(N^{1-\delta}(\log N)^{4-2\delta}) & ~ \mbox{ if } ~ \delta < 1,\\
N\frac{\zeta(s+1)\zeta(t+1)}{\zeta(s+t+2)} \sigma_{-(s+t+1)}(h) +
O(\log^3 N) & ~ \mbox{ if } ~ \delta = 1,\\
N\frac{\zeta(s+1)\zeta(t+1)}{\zeta(s+t+2)} \sigma_{-(s+t+1)}(h) +
O(1) & ~ \mbox{ if } ~ \delta > 1.
 \end{cases}
 $$ 
 where $\sigma_s(n):=\sum_{d|n}d^s$.
\end{cor}

This result has been stated in \cite{AEI} only in the asymptotic form.

\begin{cor}\label{c2}
 Let
 $$
 \phi_s(n):= n^s \prod\limits_{\substack{p|n\\ p~prime}}(1-p^{-s}).
 $$
Then for $s,t > 0$ and a positive integer $h$, we have
 $$
 \sum_{n \le N} \frac{\phi_s(n)}{n^s} \frac{\phi_t(n+h)}{(n+h)^t}
 =\begin{cases}
N\Delta(h) +
O(N^{1-\delta}(\log N)^{4-2\delta}) & ~ \mbox{ if } ~ \delta < 1,\\
N\Delta(h) +
O(\log^3 N) & ~ \mbox{ if } ~ \delta = 1,\\
N\Delta(h) +
O(1) & ~ \mbox{ if } ~ \delta > 1.
 \end{cases}
 $$
 where
 $\Delta(h)=\prod_{p|h} \left[\left(1-\frac{1}{p^{s+1}}\right)
 \left(1-\frac{1}{p^{t+1}}\right)
 + \frac{p-1}{p^{s+t+2}}\right]
 \prod_{p\nmid h}\left[\left(1-\frac{1}{p^{s+1}}\right)
 \left(1-\frac{1}{p^{s+1}}\right)
 -\frac{1}{p^{s+t+2}}\right]$,
 and $\delta=\min \{s,t\}$.
\end{cor}

\section{Preliminaries}

Some arguments for our proofs are provided by \cite{MS}, which we discuss
in the following section. In addition to that we need some well-known results
which we record below.

\begin{prop}\label{phi}
 For any real number $x \ge 1$,
 $$
 \sum_{k \le x} \varphi(k)= \frac{3}{\pi^2}x^2 + O(x \log x).
 $$
\end{prop}

\begin{defn}
 \rm
 The Mertens function $M(\cdot)$ is defined for all positive integers $n$ as
 $$
 M(n) := \sum_{k\le n} \mu(k)
 $$
 where $\mu(\cdot)$ is the M\"obius function. The above definition can be extended
 to any real number $x \ge 1$ by defining,
 $$
 M(x) := \sum_{k \le x} \mu(k).
 $$ 
\end{defn}

Essentially from the error term in the prime number theorem one gets

\begin{prop}\label{mertens}
 For any real number $x \ge 1$,
 $$
 M(x)= \sum_{k \le x} \mu(k) = O\left(x e^{-c \sqrt{\log x}}\right),
 $$
 where $c$ is some positive constant.
\end{prop}

Let $d_k(n)$ be the number of ways of writing $n$ as a product of $k$ numbers.
We generally write $d(n)$ to denote $d_2(n)$. Note that
$$
d_4(n)=\sum_{\substack{a,b \\ ab=n}} d(a)d(b).
$$
The functions $d_2(\cdot),d_4(\cdot)$ appear in our proofs.
Hence, we record the following general result about the average
order of the arithmetical function $d_k(\cdot)$, 
which can be obtained by partial summation technique.

\begin{prop}\label{k-product}
For any real number $x \ge 1$,
 $$
 \sum_{n \le x}d_k(n)= \frac{x (\log x)^{k-1}}{(k-1)!}+ O\left(x (\log x)^{k-2}\right).
 $$
\end{prop}

\section{Proofs of the Theorems}

Fine-tuning arguments of \cite{MS} does not lead us to the desired results.
However, to some extent, our proofs follow the arguments presented in \cite{MS}
verbatim, but then major steps towards our desired results are taken by
more detailed treatments for certain parts of the concerned sum. 
The proofs have two aspects of improvement. The improvement towards the
error term is obtained by a finer analysis of a particular sum and then
there is an elegant treatment of another sum
which enables us to work with the weakened hypotheses.
We elaborate below. To keep our exposition self-contained we recall relevant
parts of the proof of \thmref{error3} and \thmref{error4} from \cite{MS}.

\subsection{Setting up the proofs}

Here we explain the principle which is in the proof of \thmref{error3} and
\thmref{error4}, and also highlight the improvement towards the error term.

Let $U$ be a parameter tending to infinity which is to be chosen later. One writes,
\begin{align*}
 \sum_{n \le N} f(n) g(n) &= \sum_{n \le N} \sum_{r,s} \hat f(r) \hat g(s) c_r(n)c_s(n)\\
 &= A+B, ~\mbox{ where }
\end{align*}
$$
A:= \sum_{n \le N} \sum_{\substack{r,s\\ rs \le U}} \hat f(r) \hat g(s) c_r(n)c_s(n)
~\mbox{ and }~
B:=\sum_{n \le N} \sum_{\substack{r,s\\ rs > U}} \hat f(r) \hat g(s) c_r(n)c_s(n).
$$
To treat the sum $A$, one needs the following lemma from \cite{GMP}.

\begin{lem}\label{l1} Let $h$ be a non-negative integer. Then
 $$\sum_{n \le N} c_r(n)c_s(n+h)=\delta_{r,s} N c_r(h) + O(rs \log rs),$$
 where $\delta_{\cdot,\cdot}$ denotes the Kronecker delta function.
\end{lem}

Interchanging summations and applying \lemref{l1} (for $h=0$), one gets
upon separating $r=s$ and $r \neq s$,
\begin{align*}
 A&= N \sum_{r^2 \le U} \hat f(r) \hat g(r) \varphi(r) + O(U \log U)\\
 &=C+D+O(U \log U), ~\mbox{ where }
\end{align*}
$$
C=N \sum_r \hat f(r) \hat g(r) \varphi(r) ~\mbox{ and }~
D=-N \sum_{r^2 > U} \hat f(r) \hat g(r) \varphi(r).
$$
Clearly, $C$ is the main term as per the theorem. Using the hypotheses and
knowledge about the average order of the $\varphi$ function
(see \propref{phi}), $D$ is easily estimated (by partial summation)
to be $O \left( \frac{N}{U^\delta}\right)$. Then one estimates the sum
$B$ to obtain the result. However, we do not go into the analysis of $B$
which was done in \cite{MS}, as a major improvement that we obtain here
in this article is due to an independent treatment of the sum $B$.

When $h$ is a positive integer, the sum $\sum_{n \le N} f(n) g(n+h)$
is also written as $A+B$, but here we have
$$
A:= \sum_{n \le N} \sum_{\substack{r,s\\ rs \le U}} \hat f(r) \hat g(s) c_r(n)c_s(n+h)
~\mbox{ and }~
B:=\sum_{n \le N} \sum_{\substack{r,s\\ rs > U}} \hat f(r) \hat g(s) c_r(n)c_s(n+h).
$$
In this case, it turns out that,
$A=C+D+O(U \log U)$, where $C=N \sum_r \hat f(r) \hat g(r) c_r(h)$, the main term and
$D=-N \sum_{r^2 > U} \hat f(r) \hat g(r) c_r(h)$. To estimate $D$
one needs to know about $\sum_{r \le x} c_r(h)$, which is written as follows:
$$
\sum_{r \le x} c_r(h)=\sum_{r \le x} \sum_{d|r,d|h}\mu(r/d)d
=\sum_{\substack{k,d\\dk \le x, d| h}}d \mu(k)
=\sum_{d|h}d\sum_{k \le x/d} \mu(k).
$$
The innermost sum is $M(x/d)$, where $M(\cdot)$ denotes the Mertens
function. Using estimates on Mertens function (see \propref{mertens}),
it is then obtained that $D=O\left(\frac{N \epsilon(h)}{U^{1/2+\delta}}\right)$,
for some function $\epsilon(\cdot)$  of $h$ which is bounded above by
$e^{c\sqrt{\log h}}d(h)$ for some positive constant $c$.

However one can do better with respect to the term $O(U \log U)$ above. Note
that for the sum $A$, one actually has,
$$
A= C+D+
O \left (\sum_{\substack{r,s\\ rs \le U}} \frac{1}{(rs)^{1+\delta}} rs \log rs \right )
$$
The big-$O$ term can trivially be estimated to
be $O(U \log U)$. We note that if $\delta > 1$, the sum
$\sum_{r,s} \frac{1}{(rs)^{\delta}} \log rs$
is convergent. Hence the sum in that case is $O(1)$.
For $0<\delta \le 1$, we can write,
$$
\sum_{\substack{r,s\\ rs \le U}} \frac{1}{(rs)^{\delta}} \log rs
= \sum_{t \le U} \frac{d(t) \log t}{t^\delta}
\le \log U \sum_{t \le U} \frac{d(t)}{t^\delta}.
$$
Now from \propref{k-product} (for $k=2$) we know that
$$
\sum_{t \le U} d(t)=U \log U + O(U).
$$
Hence one can use partial summation technique to estimate the above sums.
However, estimating $\sum_{t \le U} \frac{d(t) \log t}{t^\delta}$ is
more complicated than estimating $\sum_{t \le U} \frac{d(t)}{t^\delta}$.
Since, we are only interested about the order of these sums we will
work with $\sum_{t \le U} \frac{d(t)}{t^\delta}$, as in both the cases
the resulting order is the same.

Using partial summation technique one gets,
$$
\sum_{t \le U} \frac{d(t)}{t} = O(\log^2 U)
~\mbox{ and }~
\sum_{t \le U} \frac{d(t)}{t^\delta} = O(U^{1-\delta} \log U),
$$
for $\delta < 1$.
This clearly improves the exponent of $U$ and yields,
$$
A
=\begin{cases}
C + D + O(U^{1-\delta} \log^2 U) & ~ \mbox{ if } ~ \delta < 1,\\
C+D + O(\log^3 U) & ~ \mbox{ if } ~ \delta = 1,\\
C+D + O(1) & ~ \mbox{ if } ~ \delta > 1.
\end{cases}
$$

In the following two subsections we explain our approach of handling
the sum $B$ for $h=0$ and $h \neq 0$.

\subsection{Proof of \thmref{error1}}

Recall that, by adapting the proof of \thmref{error3} one can write
$\sum_{n \le N} f(n)g(n)=A+B$, where
$$
A:= \sum_{n \le N} \sum_{\substack{r,s\\ rs \le U}} \hat f(r) \hat g(s) c_r(n)c_s(n)
~\mbox{ and }~
B:=\sum_{n \le N} \sum_{\substack{r,s\\ rs > U}} \hat f(r) \hat g(s) c_r(n)c_s(n).
$$
In the previous subsection we obtained,
$$
A
=\begin{cases}
N \sum_{r \ge 1} \hat f(r) \hat g(r) \varphi(r) +
O \left( \frac{N}{U^\delta}\right) + O(U^{1-\delta} \log^2 U)
& ~ \mbox{ if } ~ \delta < 1,\\
N \sum_{r \ge 1} \hat f(r) \hat g(r) \varphi(r) +
O \left( \frac{N}{U}\right) + O(\log^3 U)
& ~ \mbox{ if } ~ \delta = 1,\\
N \sum_{r \ge 1} \hat f(r) \hat g(r) \varphi(r) +
O \left( \frac{N}{U^\delta}\right) + O(1)
& ~ \mbox{ if } ~ \delta > 1.
\end{cases}
$$

Now here is the other major step towards improving the results of \cite{MS}.
The essence of the proof lies in a careful analysis of the term $B$.
Using \eqref{crn-exp1} and the hypothesis about $\hat f(r), \hat g(r)$ we write,
\begin{align*}
B&=\sum_{n \le N} \sum_{\substack{r,s\\ rs > U}} \hat f(r) \hat g(s) c_r(n)c_s(n)\\
& \ll \sum_{\substack{r,s\\ rs > U}}\frac{1}{(rs)^{1+\delta}}
\left| \sum_{r' | r} r' \mu(r/r') \sum_{s' | s} s' \mu(s/s')
\sum_{\substack{ m \le N/r' \\ r'm \equiv 0 \bmod s'}} 1 \right|\\
& \le \sum_{\substack{r,s\\ rs > U}}\frac{1}{(rs)^{1+\delta}}
\sum_{r' | r} r' \sum_{s' | s} s'
\sum_{\substack{ m \le N/r' \\ m \equiv 0 \bmod (s'/(r',s'))}} 1\\
& \le \sum_{\substack{r,s\\ rs > U}}\frac{1}{(rs)^{1+\delta}}
\sum_{r' | r} r' \sum_{s' | s} s' \frac{N (r',s')}{r' s'}
\end{align*}

Notice that $$\sum_{r' | r} \sum_{s' | s} (r',s') \le \sum_{l |r, l|s} l d(r/l) d(s/l).$$
This is because if $(r',s')=l$, then one can write $r'=l r''$ and $s'=l s''$.
Now since $r'|r$ and $s'|s$, we get $r'' | r/l$ and $s''| s/l$. Hence the number
of choices of $r''$ and $s''$ are at most $d(r/l),d(s/l)$ respectively.
Thus we get, upon writing $r=l r_0, ~s=l s_0$ in the sum and interchanging
summations,
$$
B \ll N \sum_{l \ge 1} \frac{1}{l^{1+2\delta}}
\sum_{\substack{r_0, s_0\\ r_0 s_0 > U/l^2}}
\frac{d(r_0) d(s_0)}{(r_0s_0)^{1+\delta}}.
$$

Now we break the outermost sum into two parts, one for $l^2 > U$ and
another for $l^2 \le U$. If $l^2 > U$, the condition $r_0 s_0 > U/l^2$ is
vacuously true and then in that case the innermost sum is a convergent
series for $\delta >0$. So we deduce that,
$$
B \ll N \sum_{l \le \sqrt{U}} \frac{1}{l^{1+2\delta}}
\sum_{\substack{r_0, s_0\\ r_0 s_0 > U/l^2}}
\frac{d(r_0) d(s_0)}{(r_0 s_0)^{1+\delta}} +
N \sum_{l > \sqrt{U}} \frac{1}{l^{1+2\delta}}.
$$
Note that the second sum is $O(\frac{N}{U^\delta})$ and also that the
innermost sum in the first sum is nothing but
$$
\sum_{t > U/l^2} \frac{d_4(t)}{t^{1+\delta}}
$$
which by partial
summation (and \propref{k-product}) turns out to be
$O\left(  \frac{\log^3(U/l^2)}{(U/l^2)^\delta}\right)$.
Putting these informations together one gets
$$
B \ll N \sum_{l \le \sqrt{U}} \frac{1}{l^{1+2\delta}}
\frac{\log^3(U/l^2)}{(U/l^2)^\delta} +
\frac{N}{U^\delta} \ll  \frac{N\log^4 U}{U^\delta}.
$$
Hence we obtain that
$$
\sum_{n \le N} f(n)g(n) 
=\begin{cases}
N \sum_{r \ge 1} \hat f(r) \hat g(r) \varphi(r) +
O \left( \frac{N\log^4 U}{U^\delta} \right) + O(U^{1-\delta} \log^2 U)
& ~ \mbox{ if } ~ \delta < 1,\\
N \sum_{r \ge 1} \hat f(r) \hat g(r) \varphi(r) +
O \left( \frac{N\log^4 U}{U}\right ) + O(\log^3 U)
& ~ \mbox{ if } ~ \delta = 1,\\
N \sum_{r \ge 1} \hat f(r) \hat g(r) \varphi(r) +
O \left( \frac{N\log^4 U}{U^\delta}\right) + O(1)
& ~ \mbox{ if } ~ \delta > 1.
\end{cases}
$$
To optimize the error terms, we choose
$$
U=\begin{cases}
N \log^2 N & ~\mbox{ if }~\delta < 1,\\
N \log N & ~\mbox{ if }~\delta = 1,\\
N^{1/\delta} (\log N)^{4/\delta} & ~\mbox{ if }~\delta > 1.
\end{cases}
$$
These choices yield us,
$$
 \sum_{n \le N} f(n) g(n) =
 \begin{cases}
N \sum_{r \ge 1} \hat f(r) \hat g(r) \varphi(r) +
O(N^{1-\delta}(\log N)^{4-2\delta}) & ~ \mbox{ if } ~ \delta < 1,\\
N \sum_{r \ge 1} \hat f(r) \hat g(r) \varphi(r) +
O(\log^3 N) & ~ \mbox{ if } ~ \delta = 1,\\
N \sum_{r \ge 1} \hat f(r) \hat g(r) \varphi(r) +
O(1) & ~ \mbox{ if } ~ \delta > 1.
 \end{cases}
 $$
This concludes the proof of \thmref{error1}.

\subsection{Proof of \thmref{error2}}

We have already mentioned that we are only left to do a careful analysis of a
particular sum, namely $B$. However in this case there are certain other
difficulties, which we get around by further subdividing the sum $B$
into parts concerning `higher' and `lower' values of $s$. One also needs
Ingham's result \cite{AEI} on the binary additive divisor problem.
We elaborate below.

Keeping the earlier principle in mind, we write,
$$
\sum_{n \le N} f(n) g(n+h) = A+B, ~\mbox{ where }
$$
$$ 
A:= \sum_{n \le N} \sum_{\substack{r,s\\ rs \le U}} \hat f(r) \hat g(s) c_r(n)c_s(n+h)
~\mbox{ and }~
B:=\sum_{n \le N} \sum_{\substack{r,s\\ rs > U}} \hat f(r) \hat g(s) c_r(n)c_s(n+h).
$$

As per our derivations above, we have,
$$
A
=\begin{cases}
N \sum_{r \ge 1} \hat f(r) \hat g(r) c_r(h) +
O\left(\frac{N \epsilon(h)}{U^{1/2+\delta}}\right) + O(U^{1-\delta} \log^2 U)
& ~ \mbox{ if } ~ \delta < 1,\\
N \sum_{r \ge 1} \hat f(r) \hat g(r) c_r(h) +
O\left(\frac{N \epsilon(h)}{U^{3/2}}\right) + O(\log^3 U)
& ~ \mbox{ if } ~ \delta = 1,\\
N \sum_{r \ge 1} \hat f(r) \hat g(r) c_r(h) +
O\left(\frac{N \epsilon(h)}{U^{1/2+\delta}}\right) + O(1)
& ~ \mbox{ if } ~ \delta > 1.
\end{cases}
$$

Now for the sum $B$, we have
\begin{align*}
B &=\sum_{n \le N} \sum_{\substack{r,s\\ rs > U}} \hat f(r) \hat g(s) c_r(n)c_s(n+h)\\
& \ll \sum_{\substack{r,s\\ rs > U}}\frac{1}{(rs)^{1+\delta}}
\left| \sum_{r' | r} r' \mu(r/r') \sum_{s' | s} s' \mu(s/s')
\sum_{\substack{ m \le N/r' \\ r'm \equiv -h \bmod s'}} 1 \right|
\end{align*}
Let $(r',s')=l$. Hence if we write $r'=lr''$ and $s'=ls''$ we get
$(r'',s'')=1$, i.e. $r''$ is invertible mod $s''$. Also note that $l | h$. Thus
\begin{align*}
B & \ll \sum_{\substack{r,s\\ rs > U}}\frac{1}{(rs)^{1+\delta}}
\sum_{r' | r} r' \sum_{s' | s} s'
\sum_{\substack{ m \le N/lr'' \\ r''m \equiv (-h/l) \bmod s''}} 1\\
& = \sum_{l|h}\sum_{\substack{r,s\\ rs > U \\l|r,l|s}}\frac{1}{(rs)^{1+\delta}} 
\sum_{r'' | r/l} lr'' \sum_{s'' | s/l} ls''
\sum_{\substack{ m \le N/lr'' \\ r''m \equiv (-h/l) \bmod s''}} 1.
\end{align*}
Writing $r=l r_0$ and $s=l s_0$, we get
\begin{align*}
B & \ll \sum_{l|h} \frac{1}{l^{2\delta}}
\sum_{\substack{r_0,s_0\\ r_0s_0 > U/l^2}} \frac{1}{(r_0 s_0)^{1+\delta}}
\sum_{r'' | r_0} r'' \sum_{s'' | s_0} s''
\sum_{\substack{ m \le N/lr'' \\ r''m \equiv (-h/l) \bmod s''}} 1\\
& =E+F,
\end{align*}
where
$$
E:= \sum_{l|h} \frac{1}{l^{2\delta}}\sum_{\substack{r_0,s_0\\ r_0s_0 > U/l^2}} 
\frac{1}{(r_0 s_0)^{1+\delta}} \sum_{r'' | r_0} r''
\sum_{\substack{s'' | s_0 \\ s'' \le N/lr''}} s''
\sum_{\substack{ m \le N/lr'' \\ r''m \equiv (-h/l) \bmod s''}} 1
$$
and
$$
F:= \sum_{l|h} \frac{1}{l^{2\delta}}\sum_{\substack{r_0,s_0\\ r_0s_0 > U/l^2}} 
\frac{1}{(r_0 s_0)^{1+\delta}} \sum_{r'' | r_0} r''
\sum_{\substack{s'' | s_0 \\ s'' > N/lr''}} s''
\sum_{\substack{ m \le N/lr'' \\ r''m \equiv (-h/l) \bmod s''}} 1.
$$

Note that,
$$
\sum_{\substack{ m \le N/lr'' \\ r''m \equiv (-h/l) \bmod s''}} 1 \le \frac{N}{l r'' s''} +1.
$$
Hence if $s'' \le N/lr''$, we get
$$
\sum_{\substack{ m \le N/lr'' \\ r''m \equiv (-h/l) \bmod s''}} 1 \ll \frac{N}{l r'' s''}.
$$
This yields,
$$
E \ll N \sum_{l |h } \frac{1}{l^{1+2\delta}}
\sum_{\substack{r_0, s_0\\ r_0 s_0 > U/l^2}}
\frac{d(r_0) d(s_0)}{(r_0 s_0)^{1+\delta}}.
$$
Now the right hand side of the above expression is bounded by
$$
N \sum_{l \ge 1 } \frac{1}{l^{1+2\delta}}
\sum_{\substack{r_0, s_0\\ r_0 s_0 > U/l^2}}
\frac{d(r_0) d(s_0)}{(r_0 s_0)^{1+\delta}},
$$
which is $O\left(\frac{N\log^4 U}{U^{\delta}}\right )$
as we have seen in the proof of \thmref{error1}.
Thus we finally have,
$$
E = O\left(\frac{N\log^4 U}{U^{\delta}}\right ).
$$

To treat the sum $F$, write $r_0=r''n_r$ and $s_0=s'' n_s$. We then deduce,
\begin{align*}
F & \le \sum_{l|h} \frac{1}{l^{2\delta}}\sum_{r_0,s_0}
\frac{1}{(r_0 s_0)^{1+\delta}} \sum_{r'' | r_0} r''
\sum_{\substack{s'' | s_0 \\ s'' > N/lr''}} s''
\sum_{\substack{ m \le N/lr'' \\ r''m \equiv (-h/l) \bmod s''}} 1\\
& \le \sum_{l|h} \frac{1}{l^{2\delta}} \sum_{r''}
\frac{1}{(r'')^\delta} \left( \sum_{n_r} \frac{1}{n_r^{1+\delta}}  \right)
\sum_{s'' > N/lr''} \frac{1}{(s'')^\delta} \left( \sum_{n_s}
\frac{1}{n_s^{1+\delta}}  \right)
\sum_{\substack{ m \le N/lr'' \\ r''m \equiv (-h/l) \bmod s''}} 1.
\end{align*}
The sums involving $n_r$ and $n_s$ are convergent, and so this simplifies to
\begin{align*}
F& \ll \sum_{l|h} \frac{1}{l^{2\delta}} \sum_{r''} \frac{1}{(r'')^\delta}
\sum_{s''> N/lr''} \frac{1}{(s'')^\delta}
\sum_{\substack{ m \le N/lr'' \\ r''m \equiv (-h/l) \bmod s''}} 1\\
& \le \sum_{l|h} \frac{1}{l^{2\delta}} \sum_{r''}
\frac{1}{(r'')^\delta} \left(\frac{l r''}{N}\right)^\delta
\sum_{s''} \sum_{\substack{ m \le N/lr'' \\ r''m \equiv (-h/l) \bmod s''}} 1\\
& = \frac{1}{N^\delta} \sum_{l|h} \frac{1}{l^\delta}
\sum_{r''} \sum_{s''} \sum_{\substack{ m \le N/lr'' \\ r''m \equiv (-h/l) \bmod s''}} 1.
\end{align*}
Now the sum
$$
\sum_{r''} \sum_{s''} \sum_{\substack{ m \le N/lr'' \\ r''m \equiv (-h/l) \bmod s''}} 1
= \sum_{n \le N/l} d(n) d(n+h/l).
$$
From the solution to the binary additive divisor problem due to Ingham \cite{AEI}
we know that
$$
\sum_{n \le N/l} d(n) d(n+h/l) \sim
\frac{6}{\pi^2} \sigma_{-1}(h/l) \frac{N}{l} \log^2(N/l).
$$
Thus, $F \ll N^{1-\delta} \log^2 N \sum_{l|h}\frac{\sigma_{-1}(h/l)}{l^{1+\delta}}$.
Note that $\sigma_{-1}(n)=O(\log n)$. Hence
$$
\sum_{l|h}\frac{\sigma_{-1}(h/l)}{l^{1+\delta}} =O(\log h).
$$
Hence we obtain,
$F \ll_h N^{1-\delta} \log^2 N$ where the implied constant depends on $h$.
Putting all these together and
then for the same choice of $U$ as in the proof of \thmref{error1} we get,
$$
\sum_{n \le N} f(n)g(n+h) 
=\begin{cases}
N \sum_{r \ge 1} \hat f(r) \hat g(r) c_r(h) +
O(N^{1-\delta}(\log N)^{4-2\delta})
& ~ \mbox{ if } ~ \delta < 1,\\
N \sum_{r \ge 1} \hat f(r) \hat g(r) c_r(h) +
O(\log^3 N) & ~ \mbox{ if } ~ \delta = 1,\\
N \sum_{r \ge 1} \hat f(r) \hat g(r) c_r(h) +
O(1) & ~ \mbox{ if } ~ \delta > 1.
\end{cases}
$$
This concludes the proof of \thmref{error2}.

\section{Concluding Remarks}

Thus, our work (\corref{c1}) gives an alternate and complete derivation
(with an explicit error term) of one of Ingham's result \cite{AEI},
which was treated only for $s,t > 1/2$ in \cite{MS}.
The condition in our main theorems, namely
$$\big|\hat f(r)\big|,\big|\hat g(r)\big| \ll \frac{1}{r^{1 +\delta}}$$
for $\delta > 0$ can perhaps be relaxed even further. 
This is currently being investigated in \cite{BS2}.  
There is considerable importance in this because, 
as indicated in \cite{GMP},
Ramanujan expansions can be used to formulate the Hardy-Littlewood
conjecture on twin primes. Thus, there is profound interest in exploring
whether our conditions can be relaxed.


\begin{thebibliography}{100}

\bibitem{padma}  H.G. Gadiyar and R. Padma, Ramanujan-Fourier series, the
Wiener-Khintchine formula and the distribution of prime pairs, 
{\em Physica A., \bf 269} (1999), 503-510.

\bibitem{GMP}  
H.G. Gadiyar, M. Ram Murty and R. Padma,
{Ramanujan - Fourier series and a theorem of Ingham},
{\em Indian J. Pure Appl. Math.}, {\bf 45} (2014), no. 5, 691-706.

\bibitem{hardy}  G.H. Hardy, Note on Ramanujan's trigonometrical
function $c_q(n)$ and certain series of arithmetical functions,
{\em Proc. Cambridge Phil. Soc.,  \bf 20} (1921), 263-271.

\bibitem{AEI}
A.E. Ingham, { Some asymptotic formulae in the theory of numbers}.
{\em J. London Math. Soc.}, {\bf 2} (1927), no. 3, 202-208.

\bibitem{RM}
M. Ram Murty, { Ramanujan series for arithmetical functions},
{\em Hardy-Ramanujan J.}, {\bf 36} (2013), 21-33. 

\bibitem{MS}
M. Ram Murty and B. Saha,
{ On the error term in a Parseval type formula in the theory of
Ramanujan expansions}, {\em J. Number Theory}, {\bf 156} (2015), 125-134.

\bibitem{SR}
S. Ramanujan, { On certain trigonometrical sums and their
applications in the theory of numbers}, {\em Trans. Cambridge Philos.
Soc.}, {\bf 22} (1918), no. 13, 259-276.

\bibitem{BS}
B. Saha, { On the partial sums of arithmetical functions with
absolutely convergent Ramanujan expansions} (to appear in Proc.
Indian Acad. Sci. Math. Sci.).

\bibitem{BS2}  B. Saha, A note on arithmetical
functions with absolutely convergent Ramanujan expansions, 
in preparation.

\end{thebibliography}
\end{document}